\documentclass[a4paper,11pt]{article}

\usepackage{jmlrutils}
\usepackage{mara}
\usepackage{lipsum}
\usepackage[dvipsnames]{xcolor}
\usepackage{amssymb}
\RequirePackage{placeins}
\RequirePackage{amsmath}
\RequirePackage{amssymb}
\RequirePackage{natbib}
\RequirePackage{graphicx}
\RequirePackage{url}
\PassOptionsToPackage{x11names}{xcolor}
\RequirePackage{xcolor}
\PassOptionsToPackage{algo2e,ruled}{algorithm2e}
\RequirePackage{algorithm2e}
\RequirePackage{hyperref}
\hypersetup{colorlinks,
            linkcolor=blue,
            citecolor=blue,
            urlcolor=magenta,
            linktocpage,
            plainpages=false}
\usepackage[margin=1.5in]{geometry}

\begin{document}

\title{On the Contractivity of Stochastic Interpolation Flow}
\author{Mara Daniels\thanks{Department of Mathematics, MIT, \texttt{maradan@mit.edu}.}}

\maketitle

\begin{abstract}%
    We investigate stochastic interpolation, a recently introduced framework for high dimensional sampling which bears many similarities to diffusion modeling. Stochastic interpolation generates a data sample by first randomly initializing a particle drawn from a simple base distribution, then simulating deterministic or stochastic dynamics such that in finite time the particle's distribution converges to the target. We show that for a Gaussian base distribution and a strongly log-concave target distribution, the stochastic interpolation flow map is Lipschitz with a sharp constant which matches that of Cafarelli's theorem for optimal transport maps. We are further able to construct Lipschitz transport maps between non-Gaussian distributions, generalizing some recent constructions in the literature on transport methods for establishing functional inequalities. We discuss the practical implications of our theorem for the sampling and estimation problems required by stochastic interpolation.
\end{abstract}

\definecolor{GreyBlue1}{rgb}{0.62, 0.73, 0.85} 
\definecolor{GreyBlue2}{rgb}{0.20, 0.35, 0.60} 
\newcommand{\mx}[1]{ 
 { \color{Orchid} \textsf{MD}: #1}
 }
\renewcommand{\cov}[0]{\ensuremath{\mathrm{cov}}}

\section{Introduction}
Recently in the high-dimensional sampling literature there has been significant interest in so-called \textit{dynamical sampling algorithms}. Roughly speaking, these algorithms generate a sample by approximately simulating the dynamics of a particle whose law converges to the target distribution. These dynamics may be stochastic or deterministic and they are typically driven by a velocity field that can be estimated from samples via least squares regression. The most famous example of this methodology is \textit{diffusion modeling}, where the idealized particle trajectory is a time-reversal of the Ornstein-Uhlenbeck process \citep{song2021scorebased}, which is now well understood and which enjoys end-to-end theoretical guarantees on the estimation error of its drift function \citep{oko2023diffusion,pmlr-v247-wibisono24a, koehler2023statistical, pabbaraju2023provable} and the discretization error of various score-based dynamical samplers \citep{chen2023sampling,benton2024nearly,li2024adapting}.

Our focus in this work is on a more recent approach to dynamical sampling known as \textit{stochastic interpolation}\footnote{Multiple authors concurrently proposed this approach under the names \textit{stochastic interpolation} \citep{albergo2023building}, \textit{flow matching} \citep{lipman2023flow}, and \textit{flow rectification} \citep{liu2023flow}.}. Given a base distribution $\mu_0 \in \cP_2(\R^d)$ and a target distribution $\mu_1 \in \cP_2(\R^d)$, a stochastic interpolant from $\mu_0$ to $\mu_1$ is any member of a family of interpolations $(\mu_t)_{t \in [0, 1]}$ constructed in the following manner. 

\begin{definition}[Definition 2.1, \cite{albergo2023stochasticinterpolantsunifyingframework}]\label{def:stochastic-interp}
    {Let $X_0 \sim \mu_0$ and $X_1 \sim \mu_1$ be independent random variables}. Let $I_t(x,x') : [0, 1] \times\R^d \times \R^d \to \R^d$ be such that for each $x, x' \in \R^d$, 
    \begin{align}\label{eqn:si-constraints}
        t \mapsto I_t(\cdot,\cdot) \in C^2([0, 1], (C^2(\R^d \times \R^d))^d) \qquad I_0(x,x') = x \qquad I_1(x,x') = x'
    \end{align}
    and additionally { $|\partial_t I_t( x, x')| \leq C |x - x'|$ } for some $C<\infty$. 
    {The $I_t$-stochastic interpolation between $\mu_0$, $\mu_1$ is the trajectory of densities $(\mu_t)_{t \in [0,1]}$ given by $\mu_t = \mathsf{Law}(I_t(X_0, X_1))$}.
\end{definition}
Here we call $\mathsf{Law}(X)$ the distribution of the random variable $X$ and $\cP_2(\Omega)$ is the class of absolutely continuous densities with finite second moment on a subset $\Omega$ of a Euclidean space. For a function $\phi : \R^d \to \R$ we denote by $\nabla \phi(x)$ its gradient, $\nabla^2 \phi(x)$ its Hessian, and for $f(x) : \R^d \to \R^d$ we call $D f(x)$ its Jacobian and $\nabla \cdot f(x) = \sum_{i=1}^d \partial_{x_i} f^{(i)}(x)$ its divergence. For functions with multiple arguments, subscripts indicate placement of derivatives, as in $\nabla_x f(x, y), D_y f(x, y)$, etc. For time dependent functions, we occasionally write $f_t(x)$ in place of $f(t, x)$ and $\dot{f}_t(x)$ in place of $\partial_t f(t, x)$, and we use process notation $(f_t)_{t \in \Lambda}$ or, when the time domain $\Lambda$ is clear from context, just $(f_t)$. 

The key property of a stochastic interpolant is that $(\mu_t)_{t\in[0,1]}$ can easily be realized by a dynamical sampler through a process known colloquially as `\textit{markovianization}.' Recall that the flow map $f_t^v(x)$ corresponding to a (lipschitz in space, uniformly in time) velocity field {$v_t(x) : [0, 1] \times \R^d \to \R^d$} is the unique solution to the ODE $\partial_t f_t(x) = v_t(f_t(x))$ with initial condition $f_0(x) = x$. 

\begin{proposition}[Theorem 2.6, \cite{albergo2023stochasticinterpolantsunifyingframework}]\label{prop:stochastic-interp-continuity}
Let $f^v_t(x)$ be the flow map corresponding to the velocity field,
    \begin{align*}
        v_t(x) & \coloneqq \E[\partial_t I_t(X_0, X_1) \mid I_t(X_0, X_1) = x] 
    \end{align*}
    The $I_t$-stochastic interpolant satisfies the \textit{continuity PDE} with respect to $(v_t)_{t \in [0, 1]}$, 
    \begin{align*}
        \partial_t \mu_t + \nabla \cdot(\mu_t v_t) = 0
    \end{align*}
    and consequently $\mu_t = (f_t^v)_{\#}\mu_0$. 
\end{proposition}
\textit{Proof sketch}. For an arbitrary test function $\phi \in C^\infty_c(\R^d)$, differentiate $\mu_t(\phi) \coloneqq \E[\phi(I_t(X_0, X_1))]$ in time:
\begin{align*}
    \partial_t \mu_t(\phi) & = \E_{X_0, X_1}[\langle \nabla \phi(I_t(X_0, X_1)), \partial_t I_t(X_0, X_1) \rangle]  \\
    & = \E_{X_t \sim \mu_t} [\langle \nabla \phi(X_t), \E[\partial_t I_t(X_0, X_1) \mid I_t(X_0, X_1) = X_t]\rangle] \\ 
    & = \E_{X_t \sim \mu_t}[ \langle \nabla \phi(X_t), v_t(X_t) \rangle],
\end{align*}
which is precisely the weak form of the continuity PDE. \hfill \ensuremath{\blacksquare} \medskip

Because $I_t(X_0, X_1)$ depends on its endpoints, the trajectory $t \mapsto I_t(X_0, X_1)$ may not be a Markov process. But by replacing the random variable $\partial_t I_t(X_0, X_1)$ with its conditional average $v_t(x)$, it is possible to construct a Markov process $t \mapsto f_t^v(X_0)$ whose marginal distribution at each $t \in [0, 1]$ coincides with $\mu_t$. It is also easy to estimate $v_t(x)$ from i.i.d. samples of $\mu_0$, $\mu_1$ by minimizing an empirical approximation to the risk
\begin{align*}
    \cL[\hat{v}] =  \E_{X_0, X_1} \left[ \| \partial_t I_t(X_0, X_1) - \hat{v}(I_t(X_0, X_1)) \|^2\right]
\end{align*}
over an appropriate class of models $\hat{v} \in \cV$, typically a class of neural networks. These observations together form the basis of a pipeline for dynamical sampling.

Owing to the flexibility of this simple construction, it is possible to define stochastic interpolants in vastly more general settings than Definition \ref{def:stochastic-interp}. Some notable examples include: $X_0, X_1$ may have an arbitrary coupling \citep{pooladian2023multisample,albergo2024stochastic}, $I_t(X_0, X_1)$ may be a random process satisfying \eqref{eqn:si-constraints} almost surely \citep[Section 3.1]{albergo2023stochasticinterpolantsunifyingframework}, or $X_0, X_1$ may be sampled from a non-Euclidean space such as the probability simplex \citep{stark2024dirichlet} or a Riemannian manifold \citep{kapusniak2024metric}. This motivates the investigation of the properties of this vast new family of parametrizations of dynamical samplers. 

Orthogonal to recent work on dynamical sampling, it is also valuable to understand the relationship between stochastic interpolation flow and other techniques for \textit{transportation of measures}. Perhaps the most famous example is the theory of optimal transportation, which furnishes for any two $\mu_0, \mu_1 \in \cP_2(\R^d)$ a unique map $T : \R^d \to \R^d$ which pushes forward $T_{\#}\mu_0 = \mu_1$ and with minimum \textit{transportation cost},
\begin{align}\label{eqn:monge-transport}
    T = \arg \inf \cJ[T] \coloneqq \E_{X \sim \mu_0}[\|X - T(X)\|^2] \qquad \text{ subject to } \qquad T_{\#} \mu_0 = \mu_1. 
\end{align}
The value of this infimum is equal to the squared Wasserstein-$2$ distance $W_2^2(\mu_0, \mu_1)$. We recall two seminal results in the theory of optimal transport which will serve as points of comparison for our main result. The celebrated \textit{Benamou-Brenier theorem} provides a dynamical perspective on optimal transportation.
\begin{proposition}[\cite{Benamou2000}]\label{thm:Benamous-Brenier}
    Given two probability measures $\mu_0, \mu_1$, it holds that 
    \begin{align} \label{eqn:kinetic-energy}
    & W^2_2(\mu_0, \mu_1) = \inf_{\tilde{v}_t} \left\{\int_0^1 \left\|\tilde{v}_t\right\|_{\cL^2\left(\mu_t\right)}^2 dt  : ~\partial_t \mu_t+ \nabla \cdot \left(\mu_t \tilde{v}_t \right)=0, \mu_{t=0}=\mu_0, \mu_{t=1}=\mu_1 \right\}.
    \end{align}
    Moreover, the optimizer $( v_t)$ induces a curve $(\mu_t)$ with the following properties. 
    \begin{enumerate}
        \item Conservative dynamics: there exists $\phi \in C^1([0,1], C^1(\R^d))$ such that $v_t(x) = \nabla \phi_t(x)$. 
        \item Zero acceleration: the flow $f^v_t(x)$ satisfies $\partial_t f^v_t(x) = T(x) - x$ where $T$ is the optimal transport map from $\mu_0$ to $\mu_1$. 
    \end{enumerate}
\end{proposition}
The zero acceleration property implies that $f^v_1(x) = x + (T(x)-x) = T(x)$ is the optimal transport map. It turns out that Lipschitz constant of $T$ is closely related to distributional properties of $\mu_0$ and $\mu_1$. Intuitively, if $\mu_1$ is more concentrated than $\mu_0$, then one might expect a bound on the contractivity of $T$ to hold, and this intuition can be made quantitative in terms of uniform bounds on log-concavity of $\mu_0$, $\mu_1$. 
\begin{definition}[Log-concavity] \label{def:log-concavity}
    A density $\mu(x) = e^{-V(x)}$ is \textit{log-concave} if $V(x)$ is a convex function. If $V \in C^2(\R^d)$, we say that $\mu$ is $\kappa$-log-concave if it holds $\nabla^2 V(x) \succcurlyeq \kappa I_d$, and it is $\eta$-log-convex if $\nabla^2 V(x) \preccurlyeq \eta I_d$, uniformly in $x \in \R^d$. 
\end{definition}
\begin{proposition}[\cite{Caffarelli2000}] \label{thm:Caffarelli}
    Let $\mu_0$ be a smooth and $\eta_0$-log-convex density, let $\mu_1$ be a smooth and $\kappa_1$-log-concave density, for $0 < \kappa_1 \leq \eta_0$. Then the optimal transport map from $\mu_0$ to $\mu_1$ is $\sqrt{\eta_0 / \kappa_1}$ Lipschitz.
\end{proposition}
Cafarelli originally proved this result by studying the regularity of solutions to a nonlinear second-order PDE known as the \textit{Monge–Ampère equation}. Later, \cite{Chewi2023EntropicCaffarelli} gave a direct proof by studying an entropic regularization of \eqref{eqn:monge-transport} and expressing the Jacobian of (an entropic analogue of) the transport map as a covariance matrix, which is bounded by well-known covariance inequalities that make use of log-concavity and convexity. Our results are in much the same spirit. As we discuss further in Section \ref{sec:lipschitz-transport}, analogues of Caffarelli's theorem have also been shown for the so-called Brownian transport map \citep{Mikulincer2024} and for a construction called `Föllmer Flow' \citep{ding2023samplingfollmerflow,dai2023lipschitztransportmapsfollmer}, which turns out to be a special case of Definition \ref{def:stochastic-interp} corresponding to $I_t(X_0, X_1) = (1-t)X_0 + tX_1$.

We make some additional assumptions that reflect the typical instantiations of stochastic interpolation in practice. {First, we specialize to the class of \textit{isotropic linear} stochastic interpolants, where $I_t$ is of the form 
\begin{align*}
    I_t(X_0, X_1) = \alpha_t X_0 + \beta_t X_1 \qquad \alpha_0 = \beta_1 = 1,\ \alpha_1 = \beta_0 = 0,\ \alpha_t, \beta_t > 0 \text{ for } t \in (0,1)
\end{align*}
which is a common parametrization in practice \citep{ma2024sit,albergo2023building,klein2023equivariant,shaul2023kinetic}.
We assume that $\alpha_t, \beta_t$ are twice continuously differentiable as functions from $[0, 1]$ to $\R$, as required by Definition \ref{def:stochastic-interp}.} Our second assumption is that $\mu_0 = \cN(0, I_d)$, which is also typical in practice and which plays a similar role to that of the Gaussian distribution in diffusion modeling. In Section \ref{sec:non-gaussian-endpoints} we discuss results for non-Gaussian $\mu_0, \mu_1$. The simplest version of our main result is as follows. 
\begin{theorem}\label{thm:isotropic-semi-gaussian}
    Suppose that $\mu_0 = \cN(0, I_d)$ and that for $V_1 \in C^2(\R^d)$, {$\mu_1(x) = e^{-V_1(x)}$ is a $\kappa$-log-concave density} with $\kappa \geq 0$. Let $(\mu_t)_{t \in [0, 1]}$ be the (isotropic) stochastic interpolant with coefficients $(\alpha_t)_{t \in [0,1]}, (\beta_t)_{t \in [0, 1]}$ and let $v_t$ be its drift. Then,
    \begin{enumerate}
        \item There exists $\phi_t \in C^1([0, 1], C^1(\R^d))$ so that $v_t(x) = \nabla \phi_t(x)$. 
        \item The velocity Jacobian is a bounded symmetric matrix,
        \begin{align*}
            D v_t(x) = \nabla^2 \phi_t(x) \preccurlyeq  \frac{\kappa \alpha_t \dot{\alpha}_t + \beta_t \dot{\beta}_t}{\kappa \alpha_t^2 + \beta_t^2} I_d
        \end{align*}
        \item The flow $f^v_t(x)$ admits the bound,
        \begin{align*}
            \|D f^v_t(x)\|_{\text{op}} \leq \sqrt{\alpha_t^2 + \beta_t^2/\kappa}
        \end{align*}
    \end{enumerate}
\end{theorem}
Since $\alpha_1 = 0$, $\beta_1=1$, Theorem \ref{thm:isotropic-semi-gaussian} implies that the endpoint of \textit{any} isotropic stochastic interpolant flow map $f^v_1$ satisfies a bound that is analogous to Caffarelli's theorem. It is easy to check that the bound in Theorem \ref{thm:isotropic-semi-gaussian} is tight by considering $\mu_1 = \cN(0,\kappa^{-1} I_d)$ in which case the stochastic interpolant flow is $f^v_1(x) = \kappa^{-1/2}x$. In fact, {for any isotropic stochastic interpolation between Gaussian measures $\mu_0 = \cN(0, \Sigma_0)$, $\mu_1 = \cN(0, \Sigma_1)$ with commuting $\Sigma_0, \Sigma_1$, the flow map $f^v_1(x)$ is equal to the optimal transport map. We prove this fact in Proposition \ref{prop:gaussian-si-is-ot} by evaluating the closed form solution of Gaussian SI identified by \citet{albergo2024stochastic}, Appendix A.} 

The proof of Theorem \ref{thm:isotropic-semi-gaussian} relies on the following {definition}, which expresses {the conditional density of $X_1 \mid X_t = x$} as an exponential family with parameter $x$. 
{
\begin{definition}\label{dfn:modified-tweedies}
    Let $X_t = \alpha_t X_0 + \beta_t X_1$ where $X_0 \sim \mu_0$, $X_1 \sim \mu_1$ independently. For $t \in [0, 1)$, the measure $\pi = \mathsf{Law}(X_1, X_t)$ has density, 
    \begin{align*}
        \pi(x_1, x_t) & = \frac{1}{(2\pi \alpha_t)^{d/2}} \exp \left( - V_1(x_1) - \frac{1}{2\alpha_t^2}\|x_t - \beta_t x_1\|^2 \right) 
    \end{align*}
    and so $X_1 \mid X_t = x$ admits the regular conditional density 
    \begin{align*}
        \mu_{t, x}(x_1) \coloneqq \exp \left(  - V_1(x_1) - \frac{1}{2}\frac{\beta_t^2}{\alpha_t^2}\|x_1\|^2 + \frac{\beta_t}{\alpha_t^2} \langle x_1, x \rangle - b_t(x) \right) 
    \end{align*}
    where $b_t(x)$ is the cumulant generating function,
    \begin{align*}
        b_t(x) \coloneqq \log \int \exp \left(  - V_1(x_1) - \frac{1}{2}\frac{\beta_t^2}{\alpha_t^2}\|x_1\|^2 + \frac{\beta_t}{\alpha_t^2} \langle x_1, x \rangle \right) \, d x_1.
    \end{align*}
\end{definition}
}
With isotropic coefficients, the substitution $X_0 = \alpha_t^{-1}(X_t - \beta_t X_1)$ yields 
\begin{align}\nonumber
    v_t(x) & = \E[\dot \alpha_t X_0 + \dot \beta_t X_1 \mid X_t = x] \\
    & = {\E\left[ \frac{\dot{\alpha}_t}{\alpha_t} X_t + \left( \dot{\beta}_t - \frac{\dot{\alpha}_t \beta_t}{\alpha_t}\right) X_1 \mid X_t = x \right] } \nonumber \\
    & = \frac{\dot{\alpha}_t}{\alpha_t} x + \left( \dot{\beta}_t - \frac{\dot{\alpha}_t \beta_t}{\alpha_t} \right) \left( \frac{\alpha_t^2}{\beta_t} \right) { \nabla b_t(x) } \label{eqn:modified-tweedie}
\end{align}
from which it is clear that $v_t(x) = \nabla \phi_t(x)$ for,
\begin{align*}
    \phi_t(x) = \frac{\dot{\alpha}_t}{2\alpha_t} \|x\|^2 + \left( \dot{\beta}_t - \frac{\dot{\alpha}_t \beta_t}{\alpha_t} \right) \left( \frac{\alpha_t^2}{\beta_t} \right) { b_t(x) }.
\end{align*}
To prove Theorem \eqref{thm:isotropic-semi-gaussian}, we bound the hessian of {$b_t(x)$} and then transfer to a bound on $f^v_t(x)$ by Grönwall's lemma. 

\section{Related Work}
Before discussing the proof of Theorem \ref{thm:isotropic-semi-gaussian} and its generalizations, we consider the relationship between our result and some of the existing literature. 

\subsection{Dynamical Sampling Methods}
Given a base measure $\mu_0$, typically a Gaussian $\mu_0 = \cN(0, I_d)$, and a target measure $\mu_1$, there is wide variety of techniques developed to construct an interpolating trajectory of measures $(\mu_t)_{t \geq 0}$ along with a velocity field $(v_t)_{t \geq 0}$ which can be used to implement a dynamical sampler for the trajectory. Some methods prefer to simulate stochastic dynamics, where $(v_t)_{t \geq 0}$ is the drift of an SDE, while other methods prefer to simulate deterministic dynamics as in Proposition \ref{prop:stochastic-interp-continuity}. Some examples include diffusion models \cite{song2021scorebased}, normalizing flows \cite{Papamakarios2021normalizing}, stochastic interpolation \cite{albergo2024stochastic}, and Poisson flow models \cite{xu2022poisson}, to name a few. In the literature on entropic optimal transport, dynamical sampling methods have also been developed for the purpose of sampling $(X_0, X_1) \sim \pi^\eps$ the $\eps$-regularized entropic optimal transport coupling of $(\mu_0, \mu_1)$ \cite{kassraie2024progressive,pooladian2024pluginestimationschrodingerbridges}.

If alongside $(v_t)$ one also has access to the \textit{score} $s_t(x) \coloneqq  \nabla \log \mu_t(x)$, then it is possible to arbitrarily trade off the amount of stochasticity in the discretization of a dynamical sampler. 
\begin{proposition}
    {Suppose} $(X_t)_{t \geq 0}$ is a stochastic process such that $X_t \sim \mu_t$ at each $t \geq 0$, and $X_t$ solves the SDE  
    \begin{align*}
        d X_t = v_t(X_t) \, dt + \sqrt{2} \, dW_t 
    \end{align*}
    where $(W_t)_{t \geq 0}$ is a Wiener process on $\R^d$. Then for $\eps \in [0,1]$, we have also $X_t' \sim \mu_t$, where $(X_t')_{t \geq 0}$ solves the SDE,
    \begin{align*}
        d X_t' = [v_t(X_t') - \eps {s_t(X_t')}]\, dt + \sqrt{2 (1-\eps)} d W_t.
    \end{align*}
\end{proposition}
In the context of diffusion modeling, this observation is used by so-called {`}probability flow ODE' models to convert SDE discretizations of diffusion to ODE discretizations \citep{song2021scorebased}. As discussed in \cite{chen2023probabilityflow}, given oracle access to $v_t(x)$ it is preferrable to discretize deterministic dynamics, since sampling error due to discretization scales as $O(\sqrt{d})$ in dimension rather than $O(d)$ scaling of SDE-based discretizations. Fully deterministic dynamics are not robust to errors in estimating $v_t(x)$ or $s_t(x)$ and require a stochastic \textit{corrector step} to alleviate these issues.

For the present work, it is important to note that stochastic interpolation with $\mu_0 = \cN(0, I_d)$ bears many similarities to diffusion modeling. As highlighted by \cite{albergo2023stochasticinterpolantsunifyingframework}, the law of an Ornstein-Uhlenbeck process with initial condition $X_0 \sim \mu_1$ can be realized as a stochastic interpolant with adjusted time domain, $(\tilde{\mu}_t)_{t \in [0, \infty]}$, and with $\alpha_t = e^{-t}$, $\beta_t = \sqrt{1-e^{-2t}}$, $\lim_{t \to \infty} \tilde{\mu}_t = \cN(0, I_d)$. The benefit of the stochastic interpolant framework is that optimizing over $(\alpha_t, \beta_t)$ can in fact lead to significant performance improvements \cite{ma2024sit,shaul2023kinetic}. In this case Theorem \ref{thm:isotropic-semi-gaussian} gives a sharp characterization of the contractivity of $v_t(x)$ in terms of $\alpha_t$, $\beta_t$, which we hope will be useful for {analyzing} discretizations in practice. 

\subsection{Score estimation}\label{sec:score-estimation}
Similar to the stochastic interpolant drift, the score $s_h(x) \coloneqq \nabla \log \mu_1 * \cN(0, h) (x)$ can be estimated by minimizing an empirical approximation to the risk 
\begin{align}
    \cL_{\text{sgm}, h}[\hat{s}_h] & = \E_{X_0, X_1} \left[ \|\hat{s}(X_1 + \sqrt{h}X_0) - s_h(X_1 + \sqrt{h}X_0) \|^2 \right]  \label{eqn:sbgm-1}\\ 
    & \cong \E_{X_0, X_1} \left[ \left\| h\hat{s}(X_1 + \sqrt{h}X_0) - \sqrt{h}X_0 \right\|^2 \right] \label{eqn:sbgm-2}
\end{align}
where $\cong$ means that the \eqref{eqn:sbgm-1} shares the same minimizers as \eqref{eqn:sbgm-2}. This equivalence follows from \textit{Tweedie's identity}, $h s_h(x) = x - \E[X_1 \mid X_1 + \sqrt{h}X_0 = x]$. 

For diffusion models, it was shown in the seminal work \citep{chen2023sampling} that the risk \eqref{eqn:sbgm-1} controls the discretization error for diffusion modeling with stochastic sampling. \cite{oko2023diffusion} construct an empirical risk minimization procedure for \eqref{eqn:sbgm-2} that achieves the rate $\tilde{O}(n^{-\frac{2s}{2s + d}})$ when the density $\mu_1$ is a Besov function $ \mu_1 \in B^s_{p, q}(\Omega)$ with uniformly upper and lower bounded density on its support $\Omega = [-1, 1]^d$. For an appropriate discretization, this leads to the rate $\E[\mathrm{TV}(\mu_1, \hat{\mu}_1)] = \tilde{O}(n^{-\frac{s}{2s+d}})$ where $\hat{\mu}_1$ is the the law of output samples, which matches the minimax rate for density estimation of $\mu_1$ up to log factors. 

More recently, \cite{pmlr-v247-wibisono24a} develop an estimator for the unsmoothed score $s_0(x) = \nabla \log \mu_1(x)$ that is compatible with densities which have unbounded support. Instead they require that $\mu_1$ is $\alpha$-sub-Gaussian and that $s_0(x)$ is $L$-Lipschitz. Analyzing the risk \eqref{eqn:sbgm-1} is challenging with classical techniques and the authors make use of deep results on the estimation of Empirical Bayes denoising functions developed in \citep{saha2020nonparametric,jiang2009emp-bayes}. The equivalence between \eqref{eqn:sbgm-1} and \eqref{eqn:sbgm-2} shows the correspondence between score estimation and estimating an empirical bayes {denoising function}, and the literature suggests using a regularized estimator of the form 
\begin{align*}
    \hat{s}_h^\eps(x) \coloneqq \frac{\nabla (\hat{\mu}_1*\cN(0, h))(x)}{\max(\eps, (\hat{\mu}_1*\cN(0, h))(x))}
\end{align*}
where $\hat{\mu}_1$ is an empirical measure supported on $n \gg 1$ samples. For an appropriate choice of bandwidth $h \sim n^{-2/(d+4)}$ and thresholding $\eps \sim n^{-2}$, the estimator $\hat{s}_h^{\eps}$ approximates $s_0$ at the rate { $\E\left[\|\hat{s}_h^\eps - s_0\|^2_{\cL^2(\mu_1)}\right] \leq \tilde{O}(n^{-2/(d+4)})$ } which is minimax optimal under the assumptions on $\mu_1$. 

{In our work, \eqref{eqn:modified-tweedie}} provides an analogue of Tweedie's identity that also suggests an analogous estimator for the stochastic interpolant drift,
\begin{align*}
    \hat{v}_{t, \eps}^h(x) = \frac{ \nabla \hat{\mu}_{t}^h(x)}{\max(\eps, \hat{\mu}_t^h(x) )} \qquad \hat{\mu}_t^h = \frac{1}{n} \sum_{i=1}^n \cN( \beta_t X_{1}^{(i)}, \alpha_t^2)
\end{align*}
where $X^{(i)}_1$, $i = 1 \ldots n$ are i.i.d. samples from $\mu_1$. Analyzing the risk of this estimator is outside the scope of our work, but it is reasonable to make comparisons to the corresponding estimator for diffusion models, i.e. that the relevant parameter for uniformly estimating $v_t(x)$ is its Lipschitz constant. Applied to diffusion models,\cite{pmlr-v247-wibisono24a} are able to bound $\E \mathrm{TV}(\hat{\mu}_{1,\eps}^h, \mu_1) \leq \tilde{O}(n^{-1/(d+4)})$ where $\hat{\mu}_{1,\eps}^h$ is the output of an appropriate diffusion model, under the assumption that $(s_t)_{t \geq 0}$ are uniformly $L$-lipschitz. Theorem \ref{thm:isotropic-semi-gaussian} suggests choosing $(\alpha_t, \beta_t)$ so that $\kappa \alpha_t^2 + \beta_t^2 = \kappa^{1-t}$, which leads to the uniform bound 
\begin{align*}
    \|Dv_t(x)\|_{\text{op}} \leq \left| \frac{d}{dt} \log \left( \kappa \alpha_t^2 + \beta_t^2 \right) \right| = \kappa.
\end{align*}
and which is compatible with the endpoint constraints on $(\alpha_t), (\beta_t)$. 

The Empirical bayes techniques used by \citep{pmlr-v247-wibisono24a} depend heavily on Gaussianity of the noise, which appears to be an obstacle in extending their analysis to have improved rates if $s_0(x)$ is very smooth. As we will discuss further in Section \ref{sec:lipschitz-transport}, the bounds we show in Theorem \ref{thm:non-gaussian-endpoints} for non-Gaussian $\mu_0, \mu_1$ are worse than the corresponding bounds on the Lipschitz constant of optimal transport maps. We suspect that Theorem \ref{thm:non-gaussian-endpoints} may be sub-optimal due to crude bounds used in our proof, which we can avoid when $\mu_0$ is Gaussian. It may be interesting to understand the sharp constant on the contractivity of the stochastic interpolant flow map, given that the corresponding constant is known for optimal transport \cite{Caffarelli2000,Chewi2023EntropicCaffarelli}.

\subsection{Lipschitz Transport}\label{sec:lipschitz-transport}

Lipschitz transportation maps play an important role in the theory of functional inequalities for probability measures. This is because if $\mu_1 = T_{\#}\mu_0$ where $T$ is $L$-lipschitz, then $\mu_1$ inherits certain functional inequalities which may be satisfied by $\mu_0$ with a constant depending on $L$ but not on the ambient dimension. One can therefore reduce the problem of proving dimension independent functional inequalities for a family of measures to the problem of constructing a Lipschitz transport map between any member and a target measure $\mu_0$ satisfying the desired inequalities. 

Caffarelli's theorem establishes such a transport map between a Gaussian measure and any strongly log concave measure. It was recently explored the possibility of constructing alternative transportation maps from finite dimensional \cite{Mikulincer2023} and infinite dimensional \cite{Mikulincer2024} Gaussian measures to $\mu_1$ which is either log-concave, or which has bounded support, leading to new inequalities on the Poincaré and log-Sobolev constants of (for example) convolutions of Gaussian measures with bounded mixing weights, as well as to drastically simplified proofs of existing inequalities \citep[Table 1]{Mikulincer2024}. In \citep{dai2023lipschitztransportmapsfollmer,ding2023samplingfollmerflow}, the authors construct a deterministic version of the Brownian transport map \cite{Mikulincer2024}, and prove it is Lipschitz using techniques similar to Theorem \ref{thm:isotropic-semi-gaussian}. The novelty of Theorem \ref{thm:isotropic-semi-gaussian} is that it holds for \textit{any} isotropic linear stochastic interpolant, allowing to understand how different choices of coefficients affect the contractivity of the velocity field. 

\section{Proof of Theorem \ref{thm:isotropic-semi-gaussian}}

We give a simple proof of Theorem \ref{thm:isotropic-semi-gaussian} that is based on the classical Brascamp-Lieb inequality \citep{BobkovLedoux2000,Bakry2014}. 
\begin{theorem}[Brascamp-Lieb]
    For $V \in C^2(\R^d)$, if $\mu(x) = e^{-V(x)}$ is a strictly log-concave density, then for $f : \R^d \to \R$,
    \begin{align*}
        \var_{X \sim \mu}[f] \leq \E_{X \sim \mu}\left[ \langle \nabla f(X), (\nabla^2 V(X))^{-1} { \nabla f(X) }\rangle \right].
    \end{align*}
\end{theorem}
From Definition \ref{dfn:modified-tweedies}, the velocity Jacobian reads
\begin{align*}
    D v_t(x) & = \frac{\dot{\alpha}_t}{\alpha_t} I_d + \left( \dot{\beta}_t - \frac{\dot{\alpha}_t \beta_t}{\alpha_t} \right) \left( \frac{\alpha_t^2}{\beta_t} \right) { \nabla^2 b_t(x) }
\end{align*}
where, using the fact that {$b_t(x)$} is a cumulant generating function, 
\begin{align*}
    { \nabla^2 b_t(x) }& = \left(\frac{\beta_t^2}{\alpha_t^4}\right) \cov_{\mu_{t,x}}[\mathsf{id}, \mathsf{id}] 
\end{align*} 
for the covariance operator 
\begin{align*}
    \cov_{\mu_{t,x}}[f, g] & = \E_{\mu_{t,x}}[f(X_1) g(X_1)^T] - \E_{\mu_{t,x}}[f(X_1)] \E_{\mu_{t,x}}[g(X_1)]^T.
\end{align*}
Now fix ${w} \in \R^d$ an arbitrary vector with $\|{w}\|_2 \leq 1$. Applying Brascamp-Lieb with test function $f_{w}(x) = \langle {w}, x \rangle$ yields 
\begin{align*}
    {w}^T \cov_{\mu_{t,x}}[\mathsf{id}, \mathsf{id}] {w} = \cov_{\mu_{t,x}}[f_{w}, f_{w}] \leq \E_{\mu_{t,x}}\left[ \langle {w}, (\nabla^2_{X_1}V_t(X_1, x))^{-1} {w} \rangle \right]
\end{align*}
for the potential {$V_t(x_1, x_t) \coloneqq V_1(x_1) + \frac{1}{2}\frac{\beta_t^2}{\alpha_t^2} \|x_1\|^2 - \frac{\beta_t}{\alpha_t^2}\langle x_1, x_t \rangle$} with,
\begin{align*}
    \nabla^2_{X_1} V_t(X_1, x) \succcurlyeq \left( \kappa + \frac{\beta_t^2}{\alpha_t^2}\right) I_d,
\end{align*}
so that 
\begin{align*}
    {w}^T D v_t(x) {w}  & \leq \frac{\dot{\alpha}_t}{\alpha_t}  + \left( \dot{\beta}_t - \frac{\dot{\alpha}_t \beta_t}{\alpha_t} \right) \left(\frac{\beta_t}{\alpha_t^2} \right) \left( \kappa + \frac{\beta_t^2}{\alpha_t^2} \right)^{-1} = \frac{\kappa \dot{\alpha}_t \alpha_t + \dot{\beta}_t \beta_t}{\kappa \alpha_t^2 + \beta_t^2}. 
\end{align*}
The bound on $\|D f^v_t(x)\|_{\text{op}}$ follows from Grönwall's argument which we provide for completeness in Lemma \ref{lem:gronwall}.

\section{Non-Gaussian Endpoints}\label{sec:non-gaussian-endpoints}

Although the assumption $\mu_0 = \cN(0, I_d)$ is typical in practical applications, it is of interest to check whether the techniques in the previous section can be applied to non-Gaussian $\mu_0$. Toward this direction, we are able to adapt the approach taken in Theorem \ref{thm:isotropic-semi-gaussian} to the setting when $\mu_0, \mu_1$ have uniform upper and lower bounds on their log-hessians. Our bound applies under the following mild condition on $(\alpha_t, \beta_t)$. 
\begin{definition}[Admissible Coefficients]
    We say that $\alpha_t, \beta_t \in C^2([0,1])$ are \textit{admissible} if $\alpha_t$ is strictly decreasing and if $\alpha_t^2 + \beta_t^2 = 1$ for all $t \in [0,1]$.
\end{definition}
When $\mu_0$ is not a Gaussian, the stochastic interpolant loses Property (1) of Theorem \ref{thm:isotropic-semi-gaussian}, which forces us to use a crude (and possible sub-optimal) bound on the velocity Jacobian. 
\begin{theorem}\label{thm:non-gaussian-endpoints}
    For $V_0, V_1 \in C^2(\R^d)$, let $\mu_0(x) = e^{-V_0(x)}$, $\mu_1(x) = e^{-V_1(x)}$ be probability densities satisfying, uniformly over $x \in \R^d$, 
    \begin{align*}
        \kappa_0 I_d \preccurlyeq \nabla^2 V_0(x) \preccurlyeq \eta_0 I_d \qquad 
        \kappa_1 I_d \preccurlyeq \nabla^2 V_1(x) \preccurlyeq \eta_1 I_d \qquad \kappa_0 \geq \kappa_1 \geq 0.
    \end{align*}
    Let $(\mu_t)$ be any stochastic interpolation with admissible coefficients $(\alpha_t, \beta_t)$, and let $(v_t)$ be its velocity. Then,
    \begin{align*}
        \|D v_t(x)\|_{\text{op}} & \leq \lambda_t \coloneqq \frac{\dot \alpha_t \alpha_t \kappa_1 + \dot \beta_t \beta_t \eta_0 }{(\alpha_t^2 \kappa_1 + \beta_t^2 \eta_0)^{1/2} (\alpha_t^2 \kappa_1 + \beta_t^2 \kappa_0)^{1/2}}. 
    \end{align*}
\end{theorem}
It follows immediately that $\|D f^v_t(x)\|_{\text{op}} \leq e^{- \int_0^t \lambda_s \, ds}$. To the best of our knowledge, this is the first example of a general purpose alternative to Cafarelli's theorem for constructing Lipschitz transport maps between log-concave $\mu_0, \mu_1$ when neither are Gaussian. From the following corollary, we see that this bound is finite for any admissible coefficients, as long as $\mu_0$ is strictly log-concave.
\begin{corollary}
    In the setting of Theorem \ref{thm:non-gaussian-endpoints}, if $\kappa_0 > 0$, then
    \begin{align*}
        \| D f^v_t(x) \|_{\text{op}} & \leq \left( \frac{\eta_0}{\kappa_1} \right)^{{ \frac{1}{2} \sqrt{\eta_0 / \kappa_0}} } 
    \end{align*}
\end{corollary} 
If we assume $\eta_0 = \kappa_0$ (in other words, that $\mu_0 = \cN(0, \kappa_0^{-1}I_d)$), then we recover Theorem \ref{thm:isotropic-semi-gaussian} as a special case. 

The corollary follows by observing that if $\beta_t = \sqrt{1-\alpha_t^2}$ increases monotonically, 
\begin{align*}
    \sup_{t\in[0,1]} \frac{\alpha_t^2 \kappa_1 + \beta_t^2 \eta_0}{\alpha_t^2 \kappa_1 + \beta_t^2 \kappa_0} = \frac{\eta_0}{\kappa_0}
\end{align*}
so $\lambda_s$ has the upper bound,
\begin{align*}
    \frac{\dot{\alpha}_t \alpha_t \kappa_1 + \dot{\beta}_t \beta_t \eta_0}{(\alpha_t^2 \kappa_1 + \beta_t^2 \eta_0)^{1/2} (\alpha_t^2 \kappa_1 + \beta_t^2 \kappa_0)^{1/2}} & \leq \frac{1}{2}\left( \frac{d}{dt} \log (\kappa_1 \alpha_t^2 + \eta_0 \beta_t^2) \right) \cdot \left( \frac{\eta_0}{\kappa_0} \right)^{1/2}.
\end{align*}

We now turn to the proof of Theorem \ref{thm:non-gaussian-endpoints}. It is convenient to work with a reparametrization of $X_0, X_1$ given by,
\begin{align*}
    \begin{bmatrix}
        X_t \\ R_t
    \end{bmatrix}
    \coloneqq
    \underbrace{
    \begin{bmatrix}
        \alpha_t I_d & \beta_t I_d \\ 
        \dot \alpha_t I_d & \dot \beta_t I_d 
    \end{bmatrix}
    }_{\coloneqq J_t}
    \begin{bmatrix}
        X_0 \\ 
        X_1
    \end{bmatrix}.
\end{align*}
Let us check that for admissible coefficients, $J_t$ is invertible. The ratio $\alpha_t / \beta_t$ is strictly decreasing, so
\begin{align*}
    0 > \partial_t \ln(\alpha_t / \beta_t) = \frac{\alpha_t}{\dot \alpha_t} - \frac{\beta_t}{\dot \beta_t} \implies \alpha_t \dot \beta_t - \beta_t \dot \alpha_t \not = 0. 
\end{align*}
This reparametrization gives a convenient expression for the stochastic interpolant drift in terms of conditional averages of $R_t$. 
\begin{definition} \label{prop:density-of-v}
    In the setting of Theorem \ref{thm:non-gaussian-endpoints}, let $(X_t, R_t) = J_t\, (X_0, X_1)$ with $X_0 \sim \mu_0$, $X_1 \sim \mu_1$ independently. Then $\tilde{\pi} \coloneqq \mathsf{Law}(X_t, R_t)$ has density 
    \begin{align*}
        \tilde{\pi}(x, r) & = \exp(-\tilde{V}_t(x, r))\, | \det J_t^{-1}\, | \\ 
        \tilde{V}_t(x, r) & \coloneqq  V_0\left( \frac{\dot \beta_t x - \beta_t r}{\alpha_t \dot \beta_t - \dot \alpha_t \beta_t}\right) + V_1\left( \frac{\alpha_t r - \dot \alpha_t x}{\alpha_t \dot \beta_t - \dot \alpha_t \beta_t} \right).
    \end{align*}
    and the SI velocity field $(v_t)$ is given by $v_t(x) = \E_{\tilde{\mu}_{t, x}}[R_t]$ where $\tilde{\mu}_{t, x}(\cdot)$ is the regular conditional density, 
    \begin{align*}
        \tilde{\mu}_{t, x}(r) & \coloneqq \frac{\exp(-\tilde{V}_t(x, r))}{\int \exp(-\tilde{V}_t(x, r)) \, dx}.
    \end{align*}
\end{definition}
We can prove Theorem \ref{thm:non-gaussian-endpoints} using a similar argument to bound $Dv_t(x)$. 

%
%

\vspace{1em}\noindent\textit{Proof of Theorem \ref{thm:non-gaussian-endpoints}.} Fix $x \in \R^d$ and let $R_t \sim \tilde{\mu}_{t, x}$.
Then, using that $R_t \sim \tilde{\mu}_{t, x}$ and $-\nabla_x \tilde{V}_t(x, R_t)$ are $\tilde{\mu}_{t, x}$-integrable,
\begin{align*}
    D v_t(x) & = D_x \left \{ 
                    \frac{\int r_t e^{-\tilde{V}_t(x, r_t)} \, dr_t}{\int e^{-\tilde{V}_t(x, r_t)}\, d r_t}
                    \right \} \\ 
                    & = \frac{\int r_t (- \nabla_x \tilde{V}_t(x, r_t))^T e^{-\tilde{V}_t(x, r_t)}\, dr_t}{\int e^{-\tilde{V}_t(x, r_t)}\, d r_t} \\ 
                    & \qquad - \left(
                        \frac{\int r_t e^{-\tilde{V}_t(x, r_t)}\, d r_t}{{\int e^{-\tilde{V}_t(x, r_t)}\, d r_t}}
                    \right)
                    \left(
                        \frac{\int (- \nabla_x \tilde{V}_t(x, r_t) )^T e^{-\tilde{V}_t(x, r_t)}\, d r_t}{{\int e^{-\tilde{V}_t(x, r_t)}\, d r_t}}
                    \right) \\ 
                    & = \cov_{\tilde{\mu}_{t, x}}[\mathsf{id}, - \nabla_x \tilde{V}_t(\cdot, x)]. 
\end{align*}
Now set
\begin{align*}
    X_0' & \coloneqq \frac{\dot \beta_t x - \beta_t R_t}{\alpha_t \dot \beta_t - \dot \alpha_t \beta_t} \qquad X_1' \coloneqq \frac{\alpha_t R_t - \dot \alpha_t x}{\alpha_t \dot \beta_t - \dot \alpha_t \beta_t}
\end{align*}
and call $f_u(R_t) = \langle u, R_t \rangle$, $g_w(R_t) = - \langle w, \nabla_x \tilde{V}_t(x, R_t) \rangle$. Then we have 
\begin{align*}
    u^T Dv_t(x) w & = u^T \cov_{\tilde{\mu}_{t,x}}[\id, - \nabla_x \tilde{V}_t(x, \cdot)] w  \leq \cov_{\tilde{\mu}_{t,x}} [f_u, f_u]^{1/2}  \cov_{\tilde{\mu}_{t,x}}[g_w, g_w]^{1/2}
\end{align*}
The first covariance is bounded by, 
\begin{align*}
    \cov_{\tilde{\mu}_{t,x}}[f_u, f_u] & \leq \E_{R_t \sim {\tilde{\mu}_{t,x}}}[\langle u, (\nabla^2_{R_t} \tilde{V}_t(x, R_t))^{-1}, u \rangle] \\ 
    & = \left( \alpha_t \dot \beta_t - \dot \alpha_t \beta_t\right)^2 E_{R_t \sim \tilde{\mu}_{t,x}}[\langle u, (\beta_t^2 \nabla^2 V_0(X_0') + \alpha_t^2 \nabla^2 V_1(X_1'))^{-1} u \rangle] \\ 
    & \leq \frac{\left( \alpha_t \dot \beta_t - \dot \alpha_t \beta_t \right)^2}{\kappa_0 \beta_t^2 + \kappa_1 \alpha_t^2}.
\end{align*}
For the second covariance, observe that
\begin{align*}
    \nabla_{R_t} g_w({R_t})& = - \left[ D_{R_t} \nabla_x \tilde{V}_t(x, {R_t}) \right] w \\ 
    & = \frac{1}{\left(\alpha_t \dot{\beta}_t - \dot{\alpha}_t \beta_t\right)^2 }\left( \alpha_t \dot{\alpha}_t \nabla^2 V_1(X_1') + \beta_t \dot{\beta}_t \nabla^2 V_0(X_0') \right) \\ 
    & = \frac{\alpha_t \dot{\alpha}_t}{\left(\alpha_t \dot{\beta}_t - \dot{\alpha}_t \beta_t\right)^2 } \left( \nabla^2 V_1(X_1') - \nabla^2 V_0(X_0') \right)  
\end{align*}
where we used that $\alpha_t^2 + \beta_t^2 = 0$ which implies $\alpha_t \dot{\alpha}_t = - \beta_t \dot{\beta}_t$. Plugging in,
\begin{align*}
    \cov_{ \tilde{\mu}_{t,x} }[g_w, g_w] & \leq \E_{{R_t}} [\langle \nabla_{{R_t}} g_w({R_t}), (\nabla_{{R_t}}^2 \tilde{V}_t(x, {R_t}))^{-1} \nabla_{{R_t}} g_w({R_t}) \rangle] \\ 
    & \leq \frac{(\alpha_t \dot{\alpha}_t)^2}{\left(\alpha_t \dot{\beta}_t - \dot{\alpha}_t \beta_t\right)^2 } \E_{{R_t}}\left[ 
        \langle v, \Gamma (X_0', X_1') v \rangle 
    \right] 
\end{align*}
where we have set $\Gamma(X_0', X_1') = \Gamma_{I} \Gamma_{II}^{-1} \Gamma_{I}$, 
\begin{align*}
    \Gamma_{I} =  
     \nabla^2 V_0(X_0') - \nabla^2 V_1(X_1') 
    \qquad 
    \Gamma_{II} =
    \alpha_t^2 \nabla^2 V_1(X_1') + \beta_t^2 \nabla^2 V_0(X_0').
\end{align*}
By Lemma \ref{lem:weird-lemma}, using that $\nabla^2 V_1(X_1') \succcurlyeq \kappa_1 I$, 
\begin{align*}
     \Gamma_{I} \Gamma_{II}^{-1} \Gamma_{I} & \preccurlyeq (\nabla^2 V_0(X_0') - \kappa_1 I) (\alpha_t^2  \kappa_1 I + \beta_t^2 \nabla^2 V_0(X_0') )^{-1}(\nabla^2 V_0(X_0') - \kappa_1 I)      
\end{align*}
Now we will use the fact that the scalar function $s \mapsto (s - \kappa)^2/(\lambda s + \kappa)$ is increasing for all $\kappa, \lambda \geq 0$ and $s \geq \kappa$. Applying this to each element of the spectrum of $\nabla^2 V_0(X_0')$, and using that $\kappa_0 \geq \kappa_1$, we get
\begin{align*}
    \Gamma(X_0', X_1') & \preccurlyeq \frac{(\eta_0 - \kappa_1)^2}{\kappa_1 \alpha_t^2 + \eta_0 \beta_t^2}\, I_d
\end{align*}
and using again the fact that $\alpha_t \dot \alpha_t = - \beta_t \dot \beta_t$,
\begin{align*}
    \cov_{ \tilde{\mu}_{t,x} }[g_w, g_w] & \leq \frac{1}{\left(\alpha_t \dot{\beta}_t + \dot{\alpha}_t \beta_t\right)^2 } \frac{(\alpha_t \dot \alpha_t  \kappa_1 + \beta_t \dot \beta_t \eta_0)^2}{\alpha_t^2 \kappa_1 + \beta_t^2 \eta_0}.
\end{align*}
Putting both these bounds together,
\begin{align*}
    D v_t(x) & \leq \left(
        \frac{1}{\kappa_0 \beta_t^2 + \kappa_1 \alpha_t^2}
    \right)^{1/2}
    \left(
        \frac{(\alpha_t \dot \alpha_t \kappa_1 + \beta_t \dot \beta_t \eta_0)^2}{\alpha_t^2 \kappa_1 + \beta_t^2 \eta_0}
    \right)^{1/2} = \lambda_t
\end{align*}
as in the statement.
\hfill $\blacksquare$

\section{Discussion}

Our approach to proving Theorem \ref{thm:isotropic-semi-gaussian} is flexible enough to extend to the cases when $\mu_0$ is non-Gaussian, but it is unclear if the gap in Theorem \ref{thm:non-gaussian-endpoints} can be improved or if it is intrinsic to stochastic interpolation. We nevertheless believe that it is a valuable and flexible tool for constructing Lipschitz transport maps. One especially interesting direction for future work is to try to generalize to the case where $\mu_0, \mu_1$ are supported on a connected Riemannian manifold. In this setting one must take into account the curvature of the ambient space, and constructing Lipschitz transport maps between concentrating measures is significantly more challenging, whereas the Brascamp-Lieb inequality still holds and it is easy to construct $I_t(X_0, X_1)$ using geodesics.

Orthogonal to this, it may also be interesting to identify other choices of $I_t(X_0, X_1)$ so that $v_t(x)$ is \textit{non-conservative}, i.e., it is not the gradient of a function. We believe that it may be possible to improve the Lipschitz constant using a flow that rotates. To gain intuition, consider the distributions 
\begin{align*}
    \mu_0 = \cN\left( 0, \begin{bmatrix}
        \kappa^{-1} & 0 \\ 
        0 & 1 
    \end{bmatrix}
    \right)
    \qquad 
    \mu_1 = \cN\left( 0, \begin{bmatrix}
        1 & 0 \\ 
        0 & \kappa^{-1}
    \end{bmatrix}\right)
\end{align*}
for which the both optimal transport map and the stochastic interpolation transport map are equal to 
\begin{align*}
    T(x) = \begin{bmatrix}
        \sqrt{\kappa} & 0 \\ 
        0 & \sqrt{1/\kappa} 
    \end{bmatrix}
    x
\end{align*}
whose Lipschitz constant diverges with the log-concavity parameter $\kappa \to 0$. The transport map $U(x)$ given by
\begin{align*}
    U(x) & = \begin{bmatrix}
        0 & 1 \\ 
        1 & 0
    \end{bmatrix}
\end{align*}
has Lipschitz constant $\|U\|_{\text{op}}=1$, but it is easy to see that $U(x)$ cannot be the gradient of any function $\phi : \R^d \to \R$. Evidently, it is possible to construct transport maps which have better contractivity properties by aligning the (ordered) principal components of $\mu_0$ with those of $\mu_1$. 

As we have shown, taking $\mu_0 \not = \cN(0, I_d)$ is one way to induce non-conservative dynamics, but it does not necessarily improve the contractivity of $f^v_t$. We hope that by tailoring $I_t(\cdot, \cdot)$ to $\mu_0, \mu_1$ and incorporating \textit{rotational alignment} of the two distributions, it will be possible to generically construct rotationally aligned transport maps that have \textit{stronger} contractivity properties than that of optimal transport. 

\vspace{3em}

\textit{The author wishes to thank Shrey Aryan, Tudor Manole, Philippe Rigollet, and Yair Shenfeld for many helpful discussions about this work. This material is based upon work supported by the U.S. Department of
Energy, Office of Science, Office of Advanced Scientific Computing Research, Department of
Energy Computational Science Graduate Fellowship under Award Number(s) DE-SC0023112.}

\newpage 
\bibliography{bibliography}

\newpage
\appendix

\section{Additional Lemmas}

\begin{lemma}\label{lem:gronwall}
    Let $\lambda_t : [0, 1] \to \R$ be such that 
    \begin{align*}
        \|{Dv_t(x)}\|_{\text{op}} \leq \lambda_t
    \end{align*}
    for any ${x} \in \R^d$. The flow map Jacobian satisfies the bound 
    \begin{align*}
        \|D f^v_t\|_{\text{op}} \leq e^{\int_0^t \lambda_s \, ds}. 
    \end{align*}
\end{lemma}
\begin{proof}
    Define {$\alpha_w(t) = Df^v_t(x)w$}. We have the differential inequality 
    \begin{align*}
        \partial_t \|\alpha_w(t)\|_2 & = \frac{1}{\|\alpha_w(t)\|_2} \alpha_w(t)^T \partial_t \alpha_w(t) \\
        & = \frac{1}{\|\alpha_w(t)\|_2} w^T (Df^v_t(x))^T D v_t(f^v_t(x)) (D f^v_t(x)) w \\ 
        & \leq {\lambda_t \frac{1}{\|\alpha_w(t)\|_2} w^T (Df^v_t(x))^T (Df^v_t(x)) w }\\
        & = \lambda_t \|\alpha_w(t)\|_2. 
    \end{align*}
    Since $\|\alpha_w(0)\| = \|w\|_2$ it follows that $\|\alpha_w(t)\|_2 \leq \|w\|_2\exp\left(\int_0^t \lambda_s \, ds\right)$.
\end{proof}

\begin{lemma}\label{lem:weird-lemma}
    For symmetric positive semidefinite matrices $A \preccurlyeq C$, $B \preccurlyeq D$, $A B A \preccurlyeq CDC$.
\end{lemma}
\begin{proof}
    For $R = C-A$, $S=D-B$, 
    \begin{align*}
    CDC & = (A + R)(B + S)(A + R) \\ 
        & = ABA + (ABR + RBA) + (RSA + ASR) + RBR + RSR + ASA \\ 
        & \succcurlyeq ABA. 
    \end{align*}
    since each term of the second line is a symmetric PSD matrix.
\end{proof}

\section{Stochastic interpolation and optimal transport for Gaussian endpoints}
    We will show that, if $\mu_0$ and $\mu_1$ are Gaussian measures with commuting covariance matrices, then the isotropic stochastic interpolant yields a flow map that coincides with the optimal transport map between the two Gaussians. Recall that the optimal transport map is given by 
    \begin{align*}
        T(x) & = m_1 + \Sigma_1^{1/2} \Sigma_0^{-1/2}(x - m_0).
    \end{align*}
    The stochastic interpolant flow map can also be calculated in closed form. Interestingly, it has the property that $f^v_1(x) = T(x)$, regardless of the choice of coefficients $\alpha_t, \beta_t$, which mirrors the behavior of Theorem \ref{thm:isotropic-semi-gaussian}.
    \begin{proposition}\label{prop:gaussian-si-is-ot}
        Suppose that $\mu_0 = \cN(m_0, \Sigma_0)$ and $\mu_1 = \cN(m_1, \Sigma_1)$ where $\Sigma_0, \Sigma_1$ are commuting positive definite matrices. Let $(\mu_t)_{t \in [0, 1]}$ be the isotropic stochastic interpolant with coefficients $(\alpha_t)_{t \in [0, 1]}$, $(\beta_t)_{t \in [0, 1]}$, and let $v_t, f^v_t$ be the drift and flow respectively. Then, 
        \begin{align*}
            v_t(x) & = \dot m_t + \frac{1}{2} \dot \Sigma_t \Sigma_t^{-1} (x - m_t) \\ 
            f_t^v(x) & = m_t + \Sigma_t^{1/2} \Sigma_0^{-1/2}(x - m_t)
        \end{align*}
        where $m_t = \alpha_t m_0 + \beta_t m_1$ and $\Sigma_t = \alpha_t^2 \Sigma_0 + \beta_t^2 \Sigma_1$. 
    \end{proposition}
    \begin{proof}
        In this setting, the stochastic interpolant drift $v_t(x)$ is computed explicitly in \citep[Equation (A.7)]{albergo2023stochasticinterpolantsunifyingframework}. Hence,
        \begin{align*}
            \partial_t (f^v_t(x) - m_t) & = v_t(f^v_t(x)) - \dot m_t = \frac{1}{2} \dot \Sigma_t \Sigma_t^{-1} (f^v_t(x) - m_t).
        \end{align*}
        Since $\Sigma_0, \Sigma_1$ commute, there exists an orthonormal matrix $U \in \R^{d \times d}$ so that $\Sigma_0 = U \Lambda_0 U^T$, $\Sigma_1 = U \Lambda_1 U^T$ are simultaneously diagonalized. Then $\Sigma_t = U\Lambda_t U^T$ for $\Lambda_t = \alpha_t^2 \Lambda_0 + \beta_t^2 \Lambda_1$ and, by linearity, we have 
        \begin{align*}
            \partial_t U^T(f^v_t(x) - m_t) & = \dot \Lambda_t \Lambda_t^{-1} U^T(f^v_t - m_t).
        \end{align*}
        Let $w_{i, t} = [U^T(f^v_t(x) - m_t)]_i$ be the $i$-th coordinate and let $\lambda_{i, t}$ be the $i$-th diagonal entry of $\Lambda_t$. By the previous display we have,
        \begin{align*}
            \partial_t w_{i, t} & = \frac{1}{2} \frac{\dot \lambda_{i, t}}{\lambda_{i, t}} w_{i, t} \implies w_{i, t} = \exp \left( \frac{1}{2}( \log(\lambda_{i, t}) - \log(\lambda_{i, 0}))\right)w_0 = \sqrt{\frac{\lambda_{i, t}}{\lambda_{i, 0}}} w_{i, 0}. 
        \end{align*}
        It follows that, 
        \begin{align*}
            U^T(f^v_t(x) - m_t) & = \Lambda_t^{1/2} \Lambda_0^{-1/2} U^T(f^v_t(x) - m_t)
        \end{align*}
        from which the desired representation follows by rearranging.
    \end{proof}

\end{document}